\documentclass{birkjour}
\usepackage{amsxtra}
\usepackage{amsthm,mathrsfs}
\usepackage{amsmath,amsthm,amssymb,amsrefs,bbm,color,enumitem,esint,esvect,float,graphicx,mathrsfs}
 \newtheorem{thm}{Theorem}[section]
 \newtheorem{cor}[thm]{Corollary}
 \newtheorem{lem}[thm]{Lemma}
 \newtheorem{prop}[thm]{Proposition}
 \theoremstyle{definition}
 \newtheorem{defn}[thm]{Definition}
 \theoremstyle{remark}
 \newtheorem{rem}[thm]{Remark}
 
 \numberwithin{equation}{section}

\title[RKT of Hankel Operators on Weighted Hardy spaces]{Reproducing Kernel Thesis \\ of Hankel Operators \\ on Weighted Hardy Spaces}

\author{Ana Čolović}
\address{
Department of Mathematics \& Statistics\\
Washington University - St. Louis\\
One Brookings Drive\\
St. Louis, MO USA\\
63130-4899}
\email{a.colovic@wustl.edu}


\newcommand\restr[2]{{
  \left.\kern-\nulldelimiterspace 
  #1 
  \vphantom{\Big|} 
  \right|_{#2} 
  }}

\begin{document}
\maketitle


\begin{abstract}
We study the boundedness of Hankel operators between two weighted spaces, with Muckenhoupt weights. In particular, we consider whether the Reproducing Kernel Thesis for Hankel operators generalizes to the case of two different weights. There, Hankel operators are bounded on the Hardy space if and only if they are bounded when tested on reproducing kernel functions. The supremum of testing the Hankel operator on this special class of functions is called the Garsia norm of the symbol of the Hankel operator, known to be equivalent to the BMO norm of the operator. We formulate a two-weight version of the Garsia norm and prove that testing the Hankel operator on the same reproducing kernel functions and measuring the norm in the two-weight setting is sufficient to prove that the operator is bounded. In the process,we prove that the Garsia norm is equivalent to the weighted Garsia norm, when two weights are the same, and we prove the two-weight version of the Carleson embedding theorem. 
\end{abstract}


\section{Introduction}

Hankel operators are operators whose matrices have the form $(\gamma_{j+k})_{j,k\geq0}$. In other words, their entries depend on the sum of coordinates. On the Hardy space $H^2,$ Hankel operators have another realization. There, given a function $\varphi \in L^2(\mathbb{T}),$ where $\mathbb{T}$ is a unit circle, we can define a Hankel operator $H_{\varphi}:H^2(\mathbb{T})\to H_{-}^2(\mathbb{T}),$ by $H_\varphi=\mathbb{P}_{-}(\varphi f),$ for polynomials $f$ in $H^2(\mathbb{T}).$  The function $\varphi$ is called the anti-analytic symbol of the Hankel operator. In \cite{Bon}, Bonsall first proved the following property of Hankel operators on Hardy spaces.

\begin{thm}\label{T:Bonsall}
Let $H_{\varphi}: H^2 \to H^2_{-}$ be  a  Hankel  operator, such that 
\begin{align*}
\sup_{z\in \mathbb{D}}\|H_{\varphi} k_z\|^2_{L^2(m)}&=\sup_{z\in \mathbb{D}}\int_{\mathbb{T}}|\varphi-\varphi(z)|^2P_z\:dm \\ &=\sup_{z\in\mathbb{D}}(|\varphi|^2(z)-|\varphi(z)|^2):= \|\varphi\|_{\mathcal{G}} < \infty.
\end{align*}
Then $H_{\varphi} $ is bounded. 
\end{thm}

Here, for $z\in \mathbb{D},\zeta \in \mathbb{T}$ $P_z=\frac{1-|z|^2}{|1-z\bar{\zeta}|^2},$ is the Poisson kernel, and for a function $\varphi \in L^2(\mathbb{T}),$ $\varphi(z)=\int_{\mathbb{T}}\varphi P_z \:dm,$ is the harmonic extension of $\varphi.$ Above, $|\varphi|^2(z)$ is the harmonic extension of $\restr{|\varphi|^2}{\mathbb{T}}$ evaluated at $z\in\mathbb{D}$ and $|\varphi(z)|^2$ is the square of the absolute value of the harmonic extension of $\varphi$ evaluated at $z\in\mathbb{D}.$ For $z\in \mathbb{D},$ $k_z$ belong to a special class of functions called normalized reproducing kernel functions.  We define them in the next section. 

We can understand Bonsall's theorem in two ways. 
The quantity on the right is known as the Garsia norm of a function $\varphi.$ The theorem then states that if the anti-analytic symbol of the Hankel operator is bounded in the Garsia norm, then the operator is bounded. Boundedness of the Garsia norm is equivalent to boundendess of the function in the $BMO$ space, an important space in harmonic analysis. Thus,  Bonsall's theorem is equivalent to saying that a sufficient condition for a Hankel operator to be bounded is for its anti-analytic symbol to belong to $BMO$. 

The theorem also gives a set of functions, namely reproducing kernels, such that if the Hankel operator is bounded on this set, the operator is bounded in general. The results of this type are known as a Reproducing Kernel Thesis for a given operator. As Bonsall's theorem shows, the Reproducing Kernel Thesis holds for Hankel operators on Hardy spaces. In this paper, we investigate whether the reproducing kernel thesis holds in  the case of Hankel operators between two weighted Hardy spaces, with weights in $A_p$, for $p=2.$

Here, $A_p$ weights, for $1<p<\infty$,  are exactly the weights that allow us to extend results of harmonic analysis from $L^p(\mathbb{T})$ to $L^p(\mathbb{T},w \: dm).$  We define them as those non-negative measurable functions $w$, such that if  $q$ is a  Hölder conjugate of $p$, 
$$[w]_{A_p}=\sup_{I\subset \mathbb{T}} \left(\frac{1}{|I|} \int_{I} w \:dm\right)\left(\frac{1}{|I|}\int_{I}w^{-\frac{q}{p}}\:dm\right)^{\frac{p}{q}}$$
is finite. Here the supremum is taken over all subarcs of $\mathbb{T}.$

In complex function theory, $A_p$ weights characterize those space of functions $L^p(\mathbb{T}, w\: dm)$ where trigonometric analytic polynomials are dense (see Nikolskii \cite{Nik 2002}).
In the weighted spaces, with weights in $A_p$, a weighted version of Green's theorem holds. So, these are the weights that preserve many of the important properties of $L^p(\mathbb{T})$ spaces. 
Given the ambient properties of weighted spaces, we wish to understand the behaviour of Hankel operators from a weighted Hardy space, and generalize the Reproducing Kernel Thesis to weighted settings.

In particular, we wish to understand whether the Reproducing Thesis holds when Hankel operators act between two weighted spaces with different weight functions. 

We can understand this question through the connection between Hankel operators and the  Hilbert transform.  In the unit disk, the Hilbert transform is defined as the operator $H:L^2(\mathbb{T})\to L^2(\mathbb{T})$ given by $H(u)=(\mathbb{P}_{+}u-\mathbb{P}_{-}u)-\hat{u}(0).$ The commutator of the Hilbert transform and a function $\varphi \in L^2(\mathbb{T})$ is the operator $[\varphi, H]f= \varphi H f - H (\varphi f),$ for functions $f\in L^2(\mathbb{T}).$ Using the definition of the Hilbert transform, we have that $[\varphi, H] = \mathbb{P}_{-}M_{\varphi}\mathbb{P}_{+}-\mathbb{P}_{+}M_{\varphi}\mathbb{P}_{-},$ where $M_{\varphi}$ denotes point-wise multiplication by the function $\varphi. $ We further observe that $[\varphi, H]f = H_{\varphi}(f_+)-H^*_{\overline{\varphi}}(f_{-},)$ where $f_+=\mathbb{P}_{+}f$ and $f_{-}=\mathbb{P}_{-}f.$ So the commutator decomposes into two orthogonal operators, a Hankel operator and an adjoint of a Hankel operator.  As we can see from Bonsall's theorem, the commutator of the Hilbert transform and the the multiplication operator by the function $\varphi$ is bounded if and only if $\varphi$ belongs to $BMO$. 

In 1985, Bloom \cite{Blo} considered the boundedness of commutators of the Hilbert transform in the two-weight setting. He considered the real-valued equivalent of the Hilbert transform, and showed that the boundedness of the commutator of the Hilbert transform on $L^p(\mathbb{R})$ is determined by the function belonging to the appropriate weighted $BMO$ class of functions. We state the result for $p=2.$
 \begin{thm} [Bloom, 1985]
    Let $\mu$ and $\lambda$ in $A_2$, and  $\nu=\mu^{\frac{1}{2}}\:\lambda^{-\frac{1}{2}}$. Then the commutator of the Hilbert transform and $M_{b}$, $[b,H](f)=bH(f)-Hb(f)$ is bounded from $L^2(\mu)$ to $L^2(\lambda)$ if and only if the function $b$ satisfies the condition
    $$\|b\|_{BMO(\nu)}:=\sup_{Q} \left(\frac{\int_{Q}|b-\langle b \rangle_{Q}|\:dx }{\int_{Q}\:\nu \: dx}\right)< \infty,$$
    \end{thm}
    
Above, the supremum is taken over all dyadic intervals in $\mathbb{R}.$

The connection between Hankel operators and the commutators of the Hilbert transform suggests that we can expect similar results to hold in the two-weight setting. So following Bloom's result, we expect Hankel operators to be bounded if the symbol has a bounded weighted $BMO$ norm. We prove the following result:

\begin{thm}[Two-weight Reproducing Kernel Thesis for Hankel operators] \label{T:main} 

Let $\mu, \lambda  \in A_2$, $\varphi \in H^2_{-}(\mathbb{T}) \cap L^2(\mathbb{T},\lambda \: dm)$ and $H_{\varphi}$ be  an operator defined by 
$$H_{\varphi} f = \mathbb{P}_{-}(\varphi f), \qquad f \text{ a polynomial}.$$  Then,
\begin{itemize} 
    \item [i)] \begin{align*}\sup_{z \in \mathbb{D}}\frac{\|H_{\varphi} K_z\|_{L^2(\lambda)}}{\|K_z\|_{L^2(\mu)}} &= \left(\sup_{z\in\mathbb{D}}\frac{1}{\mu(z)}\int_{\mathbb{T}}|\varphi-\varphi(z)|^2 P_z \: \lambda \:dm\right)^{\frac{1}{2}} \\  &:= \|\varphi\|_{\mathcal{G},\lambda \mu} < \infty \end{align*}
    \item[ii)] \begin{align*}
    \sup_{z \in \mathbb{D}}\frac{\|H_{\varphi} K_z\|_{L^2(\mu^{-1})}}{\|K_z\|_{L^2(\lambda^{-1})}}&= \left(\sup_{z\in\mathbb{D}}\frac{1}{\lambda^{-1}(z)}\int_{\mathbb{T}}|\varphi-\varphi(z)|^2 P_z \: \mu^{-1} \:dm\right)^{\frac{1}{2}} \\ &:= \|\varphi\|_{\mathcal{G},\mu^{-1} \lambda^{-1}}< \infty\end{align*}
\end{itemize}
if and only if $H_\varphi$ is bounded as an operator between $H^2(\mu)$ and $\overline{H^2_{-}\cap L^2(\lambda})$ and as an operator between $H^2(\lambda^{-1})$ and $\overline{H^2_{-}\cap L^2(\mu^{-1})}$
\end{thm}
\par
Here $H^2(\mu)$ is defined as the closure of polynomials in the $L^2(\mathbb{T},\mu)$ norm, and the closure of $H^2_{-} \cap L^2(\lambda)$ is taken in the $L^2(\lambda)$ norm. Similar definitions hold for the $\lambda^{-1}$ and $\mu^{-1}$ case. 
Note that $\mathbb{P}_{+}:L^2(\mathbb{T},\lambda)\to L^2(\mathbb{T},\lambda)$ is bounded if and only if $\lambda$ belongs to $A_2$ (see \cite{Nik 2002}). Using this fact, and the assumption that $\varphi \in L^2(\mathbb{T},\lambda),$ it follows that the Hankel operator $H_{\varphi}$ above, is well-defined on polynomials. The definition of the operator then extends to all continuous functions on the disc, by uniform approximation. For the forward direction, we prove that the operator is bounded on the polynomials, and then equate boundedness of the operator with boundedness on a dense class. As we will see, the functions $K_z$ are continuous and belong to $H^2(\mu)$ and $H^2(\lambda).$ So, the reverse direction follows by testing on $K_z,$ for each $z\in \mathbb{D}.$

 Above, $K_z$ are reproducing kernel functions, for $z\in\mathbb{D}$. We note that for $z\in\mathbb{D}$, $\varphi K_{z}=\mathbb{P}_{-}(\varphi K_{z})+\mathbb{P}_+(\varphi K_{z}).$
 According to calculations in \cite{Tre}, $\mathbb{P}_+(\varphi K_z)=\varphi(z)K_z$, and so $H_{\varphi} K_{z}=\mathbb{P}_{-} (\varphi K_{z})= (\varphi-\varphi(z))K_{z}. $ Thus, 
$$\sup_{z\in \mathbb{T}}\frac{\|H_{\varphi} K_{z}\|_{L^2(\lambda)}^2}{\|K_{z}\|_{L^2(\mu)}^2}=\sup_{z\in\mathbb{D}}\frac{\int_{\mathbb{T}}|\varphi-\varphi(z)|^2 |K_{z}|^2 \lambda \:dm}{\int_{\mathbb{T}}|K_{z}|^2 \mu \:dm}.$$
Note that for $z,w\in\mathbb{D}$, $P_{z}(w)=\frac{1-|z|^2}{|1-z\bar{w}|^2}=(1-|z|^2)|K_{z}|^2.$ So, multiplying the top and the bottom of the above equations by $1-|z|^2,$ for each $z\in\mathbb{D}$, we get that
\begin{equation}\label{Gamma-Garsia-weighted}
\sup_{z\in\mathbb{D}}\frac{\|H_{\varphi} K_z\|^2_{L^2(\lambda)}}{\|K_z\|^2_{L^2(\mu)}} =  \sup_{z\in\mathbb{D}}\frac{1}{\mu(z)}\int_{\mathbb{T}}|\varphi-\varphi(z)|^2P_z \lambda \:dm.
\end{equation}
We also note that:
\begin{equation}\label{Gamma-Garsia-conjugate-weighted}
\sup_{z\in\mathbb{D}}\frac{\|H_{\varphi} K_z\|^2_{L^2(\mu^{-1})}}{\|K_z\|^2_{L^2(\lambda^{-1})}} =  \sup_{z\in\mathbb{D}}\frac{1}{\lambda^{-1}(z)}\int_{\mathbb{T}}|\varphi-\varphi(z)|^2P_z \mu^{-1} \:dm,
\end{equation}
for inverse weights of $\mu,$ and $\lambda.$
We say that a function has a bounded two-weight Garsia norm if both the quantity on the right hand side of \eqref{Gamma-Garsia-weighted} and \eqref{Gamma-Garsia-conjugate-weighted} are bounded.
We denote the quantity (\ref{Gamma-Garsia-weighted}) by $\|\varphi\|_{\mathcal{G}, \lambda\mu}.$ 
and the quantity (\ref{Gamma-Garsia-conjugate-weighted}) by $\|\varphi\|_{\mathcal{G}, \mu^{-1}\lambda^{-1}}.$ 
Bonsall's theorem states that Hankel operator is bounded if and only if its anti-analytic symbol has bounded Garsia norm. The theorem above states that if the anti-analytic symbol has a bounded two-weight Garsia norm, then the associated Hankel operator is bounded. 

When proving the above theorem, our approaches differ in the case when $\mu$ and $\lambda$ are the same weight, and in the general case. We present these cases separately.

In the case that $\mu=\lambda,$ we prove that quantities (\ref{Gamma-Garsia-weighted}), and (\ref{Gamma-Garsia-conjugate-weighted}) are equivalent to the unweighted Garsia norm, which is known to be equivalent to the $BMO$ norm. In this case, the main theorem states that boundedness of Hankel operators in the weighted case can be inferred from the boundedness of its anti-analytic symbol in $BMO$, just as in the unweighted case.
In \cite{Tre}, Treil showed that Bonsall's theorem holds using Green's theorem and a closely related lemma by Uchiyama. To prove the above theorem in the one-weight setting, we follow Treil's approach.  
In this case, an integral part of aligning the weighted and the unweighted setting, is proving the equivalence of one-weight Garsia norm and the Garsia norm of a function, mentioned above. This allows us to prove that a certain embedding operator is bounded, and as a consequence prove the main theorem.

For a function $\varphi$ the Garsia norm is defined to be
 $$\|\varphi\|_{\mathcal{G}}=\sup_{z\in \mathbb{D}}\left(\int_{\mathbb{T}}|\varphi-\varphi(z)|^2P_z\:dm\right)^{\frac{1}{2}}.$$ We define $$\|\varphi\|_{\mathcal{G},w}:=\sup_{z\in\mathbb{D}}\left(\frac{1}{w(z)}\int_{\mathbb{T}}|\varphi-\varphi(z)|^2P_zw\:dm\right)^{\frac{1}{2}}$$ to be its one-weight Garsia norm.
 
 We state the theorem that concerns the equivalence of one-weight and the unweighted Garsia norms of the function. In proving it, we use ideas of \cite{Nik 2002}.
\begin{thm}\label{T:Garsia-weighted} 
Suppose $w\in A_2.$ For a function $\varphi \in L^2(\mathbb{T}),$ the weighted $BMO$ and the weighted Garsia norms are equivalent, i.e.
$$ \|\varphi\|_{\mathcal{G},w} \approx \|\varphi\|_{\mathcal{G}}$$
where the implied constants depend only on $w.$
\end{thm}

In the two-weight case, as above, we need to prove that an appropriate embedding operator is bounded. The following theorem acts as a bridge between our assumptions on the anti-analytic symbol and the boundedness of the appropriate embedding operator. It is a generalization of the Carleson embedding theorem to the two-weight setting. 

\begin{thm} \label{T:Two-weight-Carleson} Let $d\tau$ be a nonnegative measure in $\mathbb{D}.$ Then the following assertions are equivalent.

\begin{itemize}
    \item[i)] For all $f\in L^2(\mathbb{T}, \mu),$
    $\int_{\mathbb{D}}|f|^2 \lambda \: d\tau \leq B \int_{\mathbb{T}}|f|^2 \mu \: dm$
    \item [ii)] $\sup_{z \in \mathbb{D}} \frac{1}{\mu(z)}\|P_{z}\|_{L^1(\lambda d\tau)} < C <\infty$
    \item [iii)] There is a constant $D$, such that $\frac{1}{\mu(I_{S})} (\lambda d\tau)(S) < D$  for all sectors $S=\{re^{i\theta}: 1-h \leq r < 1, |\theta-\theta_0|\leq h\},$ where $I_{S}$ is an interval of length $h$ centered at $e^{i\theta_{0}}.$
\end{itemize}
If any of the assertions hold true, we say that the measure $d\tau$ is $\lambda-\mu$ Carleson.
\end{thm}

Throughout the paper, as well as in the above theorem, we use the standard notation for equivalence of two quantities. Namely, for quantities $A$ and $B$, we use $A\lesssim B$ to indicate that $A \leq CB$, for some constant $C$. We use $A \approx B,$ to mean that $A \lesssim B,$ and $B \lesssim A.$

\par We organize the paper in the following way. In the next section, we give the background necessary for proving the main theorem. Here we also prove the equivalence of Garsia and the one-weight Garsia norms. In Section 3, we prove the main theorem in the one-weight setting. In Section 4, we prove the two-weight Carleson embedding theorem, and as its consequence we prove the main theorem in full generality. Finally, in the last section, we prove  a consequence of the two-weight Carleson embedding theorem, namely boundedness of a certain class of integral operators between weighted spaces, with two weights. 

%
\section {Preliminaries and Theorem \ref{T:Garsia-weighted}}

\subsection{$A_p$ weights}

Since we aim to generalize Bonsall's theorem to weighted settings, we begin by introducing $A_p$ weights and discussing their relevant properties. What follows are standard results, so for each of them, we  provide either a brief proof or a reference for interested readers. 

We begin with a definition of $A_p$ weights. 
\begin{defn} We say that a non-negative measurable function $w$, belongs to $A_p,$ for $1<p<\infty,$ and $q$ its Hölder conjugate if 
$$[w]_{A_p}=\sup_{I\subset \mathbb{T}} \left(\frac{1}{|I|} \int_{I} w \:dm\right)\left(\frac{1}{|I|}\int_{I}w^{-\frac{q}{p}}\:dm\right)^{\frac{p}{q}}$$
is finite. Here the supremum is taken over all subarcs of $\mathbb{T}.$ We define $A_\infty$ to be the union of $A_p$ for $1<p<\infty.$
\end{defn}

We suppose $w\in A_2.$ Throughout the paper we use the following notation.  For $I\subset \mathbb{T},$

$$ w(I) = \int_I w \:dm, \qquad w_I = \frac{w(I)}{|I|}.$$

We will make frequent use of the following lemma which establishes the equivalence of $[w]_{A_2}$ and a similar Poisson condition of a weight. Its proof can be found in \cite{Nik 2002}. 

\begin{lem}\label{L:Poisson-Ap} 
Suppose $w \geq 0,$ $w,\frac{1}{w} \in L^1(\mathbb{T}).$ Then the following statements are equivalent.
\begin{itemize}
\item[(i)] $[w]_{A_2}=:\sup_{I\subset \mathbb{T}}\left(\frac{1}{|I|}\int_I w \:dm\right)\left(\frac{1}{|I|} \int_I w^{-1} \:dm\right) < \infty  $.
\item[(ii)] $PA_2(w)=:\sup_{z\in \mathbb{D}} w(z) \frac{1}{w}(z) < \infty $.

\end{itemize}
\end{lem}

If a weight belongs to the $A_p$ class of weights, then it also belongs to $A_{p-\epsilon}$, for some $\epsilon >0.$ We make note of this property of $A_2$ weights, found in \cite{Nik 2002}, as we use it in sections that follow.

\begin{lem}\label{L:epsilon-of-room} 
Let $w\in A_2$. Then there exists $\epsilon >0,$ such that
$w\in A_{2-\epsilon}.$
\end{lem}

$A_\infty$ weights also have the following useful property.

\begin{lem}\label{L:Inverse-cond-for-Ap} 
Suppose that $w \in A_p,$ for  $1<p<\infty.$ Then
 for every subarc $I$ of $\mathbb{T}$, and every measurable subset $E$ of $I,$
\begin{equation*}
\left(\frac{|E|}{|I|}\right)^p \leq [w]_{A_p} \frac{\int_E w \:dm}{\int_I w \:dm}. \label{A_infinity}
\end{equation*}

\end{lem}

The proof of the above lemma is identical as for the case of finite intervals in $\mathbb{R}^n$, found in \cite{GunWhe}. Note that by Lemma \ref{L:epsilon-of-room}, if $w\in A_2$ the results above holds with $\eta=2-\epsilon.$

As a consequence of Lemma \ref{L:epsilon-of-room}, we have the following doubling property of $A_2$ weights. 

\begin{lem}\label{L:inverse-doubling-w} 
Let $w\in A_2.$ Suppose $I$ and $J$ are arcs in $\mathbb{T}$, with the same center and $|J|=2|I|$. Then $\frac{1}{w(I)}\leq \frac{4[w]_{A_2}}{w(J)}$.
\end{lem}

\begin{proof}
By Lemma \ref{L:inverse-doubling-w}, since $w\in A_2$ it follows that $\left(\frac{|I|}{|J|}\right)^2\leq [w]_{A_2}\frac{\int_I w\:dm}{\int_J w\:dm}.$ Thus $\frac{1}{4}\leq [w]_{A_2}\frac{\int_I w\:dm}{\int_J w\:dm},$ and so $\frac{1}{\int_I w\:dm}\leq 4[w]_{A_2}\frac{1}{\int_J w\:dm}.$ 
\end{proof}

\subsection{$BMO$ and one-weight $BMO$ norms}

Properties of Hankel operators $H_{\varphi}$ are often tied to properties of their anti-analytic symbols, $\varphi.$ In the case of Hardy spaces,  boundedness of the operator is related to boundedness of the anti-analytic symbol in the $BMO$ norm. As Theorem \ref{T:main} will establish, the corresponding result also holds in weighted Hardy spaces. In this section, we define $BMO$ and one-weight $BMO$ spaces, and discuss some of their relevant properties. In what follows, we refer to the one-weight spaces as weighted ones, and in later sections make a distinction between two-weight spaces.

\par A function $\varphi$ is said to be of bounded mean oscillation, and belong to $BMO$ if the following quantity is finite
$$\|\varphi\|_{BMO}:=\sup_{I \subset \mathbb{T}} \frac{1}{|I|} \int_{I} |\varphi-\varphi_{I}|\:dm.$$
Here, the supremum is taken over all subarcs $I$ of $\mathbb{T},$ and $\varphi_I$ is defined as $\varphi_I=\int_I \varphi \:dm.$ 
Here, the averages are taken with respect to the indicator function kernels. It is also possible to define averages with respect to Poisson kernels. This norm is the Garsia norm, and for a function $\varphi$ it is defined as
$$\|\varphi\|_{\mathcal{G}}=\sup_{z\in \mathbb{D}}\left(\int_{\mathbb{T}}|\varphi-\varphi(z)|^2P_z\:dm\right)^{\frac{1}{2}}=\sup_{z\in\mathbb{D}}(|\varphi|^2(z)-|\varphi(z)|^2)^{\frac{1}{2}}.$$
It is a well known fact that $BMO$ and Garsia norms are equivalent. A proof of this can be found in \cite{Gar}.

In \cite{MucWhe}, Muckenhoupt and Wheeden defined the weighted $BMO$ norm, and proved its equivalence to the $BMO$ norm. 
If $w\in A_2,$ we say that a function $\varphi$ belongs to the weighted $BMO$ space, $BMO_w$  if the following quantity is finite:
$$\|\varphi\|_{BMO, w}:=\sup_{I\subset \mathbb{T}}\frac{1}{\int_I w \:dm } \int_{I} |\varphi-c_{I}|\:dm,$$
where for $I \subset \mathbb{T},$ $c_I=\frac{1}{\int_I w \:dm} \int_I \varphi \:dm.$ 

For $w\in A_2$, we define the weighted Garsia space analogously to the unweighted case. We say that a function $\varphi$ belongs to $BMO_{\mathcal{G},w} $ if  the following quantity is finite. 
$$\|\varphi\|_{\mathcal{G},w}:=\sup_{z\in \mathbb{D}}\left(\frac{1}{w(z)}\int_{\mathbb{T}}|\varphi(\zeta)-\varphi(z)|^2P_z(\zeta)w(\zeta)\:dm(\zeta)\right)^{\frac{1}{2}}.$$

Collecting the equivalences of norms mentioned above, for a function $\varphi$, we have the following set of equivalences:
\begin{align}
    \|\varphi\|_{BMO, w} &\approx \|\varphi\|_{BMO} \label{BMO-weighted-BMO} \\
    \|\varphi\|_{BMO} &\approx \|\varphi\|_{\mathcal{G}} \label{BMO-Garsia}
\end{align}

From these equivalences we can deduce that $\|\varphi\|_{BMO, w} \approx \|\varphi\|_{\mathcal{G}},$ i.e. weighted BMO and the Garsia norms are equivalent.  

The main step in proving Theorem \ref{T:main} will the proving the equivalence of Garsia and the weighted Garsia norms. In light of norm equivalences mentioned so far, we prove this result as a simple corollary of the following lemma.

\begin{lem}\label{L:weighted-BMO-weighted-Garsia} 
Suppose $w\in A_2.$ For a function $\varphi,$ the weighted BMO and the weighted Garsia norms are equivalent, i.e.
$$ \|\varphi\|_{\mathcal{G}, w} \approx \|\varphi\|_{BMO, w},$$
where the implied constants depend only on $w.$
\end{lem}

By Remark 2.4, and equivalences \eqref{BMO-weighted-BMO} and \eqref{BMO-Garsia}, Theorem \ref{T:Garsia-weighted} follows as a corollary.

For functions in $BMO,$ we have the following consequence of John-Nirenberg inequality, found in \cite{Gar}. If $\varphi \in BMO$, then for any finite $p>1$, \begin{equation*}\|\varphi\|_{BMO} \approx \sup_{I\subset \mathbb{T}}\left(\frac{1}{|I|}\int_{I}|\varphi-\varphi_{I}|^p\:dm\right)^{\frac{1}{p}}. \end{equation*}

In \cite{Ho}, Ho established the equivalent result in $BMO_w$ spaces.

\begin{lem}\label{L:JN-weighted} 
Suppose $\varphi \in BMO_w$ for $w\in A_{\infty}.$ Then for $0< p<\infty,$ 
\begin{equation}\label{bmo-L-p}\|\varphi\|_{BMO,w} \approx \sup_{I\in \mathbb{T}} \left(\frac{1}{w(I)}\int_{I}|\varphi - c_{I}|^p\:dm\right)^{\frac{1}{p}}, 
\end{equation} 
where $c_I=\frac{1}{\int_I w \:dm} \int_I \varphi \:dm,$ for $I\subset \mathbb{T}.$
Namely, there exist constants $C'_1(w)$ and $C'_2(w)$ that only depend on $w$ such that $$C'_1(w)\|\varphi\|_{BMO,w} \leq \sup_{I\in \mathbb{T}} \left(\frac{1}{w(I)}\int_{I}|\varphi - c_{I}|^p\:dm\right)^{\frac{1}{p}}\leq C'_2(w)\|\varphi\|_{BMO,w}.$$  

\end{lem}

In \cite{Ho}, the lemma is stated as an equivalence of the $BMO$ norms and the quantity on the left of (\ref{bmo-L-p}).  Since $BMO$ and weighted $BMO$ norms are equivalent, we state it as the equivalence of the weighted $BMO$ norm and the quantity above. 

In proving Theorem \ref{T:main}, we also use the following result, due to Wang \cite{Wan}, that allows for greater freedom related to constants, when working in $BMO_w$ spaces.

\begin{lem}\label{L:JN-w-c} 
Suppose $\varphi \in BMO_w$, for $w\in A_{\infty}.$ Then, for $0< p<\infty,$
$$||\varphi||_{BMO,w} \approx \underset{I\subset \mathbb{T}}{\sup} \ \underset{c}{\inf}\left(\frac{1}{w(I)}\int_I |\varphi -c|^pw\:dm\right)^{\frac{1}{p}} $$
where the infimum above is taken over all the positive constants $c,$ and and the implied constants depend only on $w.$
\end{lem}

\begin{rem} \label{R:induction-BMO}
We observe that the difference of weighted averages can be bounded by the weighted $BMO$ norm of a function. 
Suppose that $I,J\subset \mathbb{T}$ are  such that $|J|=2|I|.$ Then the difference of averages can be estimated using Lemma \ref{L:inverse-doubling-w}. Namely, 
\begin{align}
 |c_I-c_J|& \leq \frac{1}{w(I)}\int_I|\varphi-c_J|w\:dm \notag \\
 &\leq \frac{4[w]_{A_2}}{w(J)}\int_J|\varphi-c_J|w\:dm \notag \\
 &\leq 4[w]_{A_2} \|\varphi\|_{BMO, w}. \notag
\end{align} 
As an extension of this observation, for an arc $I$ as above and $J=2^kI,$ where $k \geq 1$, we have, using induction that
\begin{align} \label{averages}
|c_I-c_J|&\leq |c_{2^kI}-c_{2^{k-1}I}|+|c_{2^{k-1}I}-c_{2^{k-2}I}|+ \dots + |c_{2I}-c_I|\notag \\ &\leq k\, 4 \, [w]_{A_2}\|\varphi\|_{BMO, w}. 
\end{align}
\end{rem}


\subsection{Green's theorem, Uchiyama's lemma}

In this section, we collect some well known results about function theory in the unit disc. These are standard results in the area, so the proofs are omitted. As before, we provide references for interested readers. 

We use $dm$ to denote the normalized arc length measure on the circle, and $dA$ the normalized area measure on the unit disk. 

The following is the standard version of Green's theorem, the proof of which can be found in Zygmund \cite{Zyg}.

\begin{lem}[Green's Theorem]\label{Green's Theorem}
Let $f \in C^2(\mathbb{D})\cap C(\overline{\mathbb{D}} ).$ Then
$$\int_{\mathbb{T}} f(\zeta)\:dm(\zeta) - f(0) = \int_{\mathbb{D}} \Delta f(z)  \ln\frac{1}{|z|}\:dA(z).$$
\end{lem}

The theorem below can be found in Garnett \cite{Gar}. It gives an equivalence of several conditions, each of which then defines a Carleson measure. 

\begin{thm}[The Carleson embedding theorem] \label{T:Carleson} 
Let $d\tau$ be a nonnegative measure in $\mathbb{D}.$ Then the following assertions are equivalent.
\begin{itemize}
    \item[i)] For all $f\in L^2(\mathbb{T}),$
    $\int_{\mathbb{D}}|f|^2 \: d\tau \leq C \int_{\mathbb{T}}|f|^2 \: dm$
    \item [ii)] $\sup_{z \in \mathbb{D}} \|P_{z}\|_{L^1(d\tau)} < \infty$
    \item [iii)] There is a constant $c$, such that $\frac{1}{h} d\tau(S) < c$  for all sectors $S=\{re^{i\theta}: 1-h \leq r < 1, |\theta-\theta_0|\leq h\}.$
\end{itemize}
\end{thm}

The following lemma is closely related to Green's theorem. Its proof can be found in \cite{Nik 1986}.

\begin{lem}[Uchiyama's lemma]\label{L:Uchiyama} 
Let $u \in C(\overline{\mathbb{D}} )$, be a non-negative subharmonic function ($\Delta u \geq 0$). Suppose $0\leq u(z) \leq 1$, for $z \in \mathbb{D}$. Then the measure $\mu(z)=\Delta u(z) \ln\frac{1}{|z|}\:dA(z)$ is a Carleson measure.
\end{lem}

The lemma above states that if  $u$ is a subharmonic function,  then for any function $f \in L^2(\mathbb{T}),$
$$\int_{\mathbb{D}}|f(z)|^2\Delta u (z) \:  \ln\frac{1}{|z|}\:dA(z)\leq e\int_{\mathbb{T}}|f(\zeta)|^2\:dm(\zeta).$$

We use the following lemma in both the one-weight and two-weight settings. 
In the one-weight setting, its corrolary allows us to 
translate the properties of a Carleson measure to  $L^2(w\:dm)$ spaces, with $w\in A_2$, and in the two-weight setting it forms an important piece of our weighted Carleson theorem. It was originally proved by Treil and Volberg and it can be found in \cite{Nik 2002}.

\begin{lem}[Treil and Volberg, 1997] \label{L:One-weight-Treil-and-Volberg} 
Let $w\geq 0, $ $w, \frac{1}{w} \in L^1(\mathbb{T}),$ and let $dv$ be a non-negative measure on $\mathbb{D}$ such that 
$$\int_{\mathbb{D}}P_z(\xi)(w(\xi))^2 \: dv(\zeta) \leq B^2\, w(z)= B^2 \int_{\mathbb{T}}P_z(\zeta)w(\zeta)dm(\zeta)$$
for every $z \in \mathbb{D},$ and for a constant $B>0.$ Then
$$\int_{\mathbb{D}}|f(\zeta)|^2\: dv (\zeta) \leq 16 B^2 \int_{\mathbb{T}}|f|^2w^{-1}\: dm,$$ for all functions $f\in L^2(\mathbb{T},w^{-1}).$
\end{lem}

The following corollary says $A_2$ weights preserve the embedding property of Carleson measures. 

\begin{cor}\label{L:weighted-carleson} 
Let $w \in A_2.$ Suppose $\mu$ is  Carleson measure on 
$\mathbb{D}.$ If $f\in L^2(\mathbb{T},\:w),$ then 
$$\int_{\mathbb{D}}|f(z)|^2w(z)d\mu(z)\leq C_\mu(w)\int_{\mathbb{T}}|f|^2w\:dm,$$
where $C_\mu(w)$ is a constant that depends on $w$ and the Carleson embedding constant $C_{\mu}.$
\end{cor}

In the proof of Theorem \ref{T:main}, we also use the following lemma, which can be found in \cite{Nik 2002}. 

\begin{lem}\label{L:Green-weighted} 
Suppose $w\in A_2.$ Then 
$$\frac{1}{c(w)^2}\int_{\mathbb{T}}|f|^2w\:dm\leq \int_{\mathbb{D}}|\nabla f(z)|^2 w(z)  \ln\frac{1}{|z|} \:dA(z) \leq c(w)^2 \int_{\mathbb{T}}|f|^2 w\:dm, $$
for every f $\in L^2(\mathbb{T}, w)$ such that $\hat{f}(0)=0$, where $c(w)$ is a constant that depends on $w.$

\end{lem}

\subsection{Proof of Theorem \ref{T:Garsia-weighted}}

We now prove Theorem \ref{T:Garsia-weighted}, establishing the equivalence of Garsia and one-weight Garsia norms. As discussed above, we do so by first proving Lemma \ref{L:weighted-BMO-weighted-Garsia}, namely the equivalence of weighted $BMO$ and weighted Garsia norms. 

We let $z_I\in \mathbb{D}$ be a point such that $1-|z|=|I|$, with $z/|z|$ the center of $I$. Note that for all $z\in \mathbb{D},$ there exists $I\subset \mathbb{T}$ such that $z=z_I$. 

We first prove that  
$$\|\varphi\|_{BMO,w} \lesssim  \|\varphi\|_{\mathcal{G},w}.$$

Assume $\varphi \in BMO_{\mathcal{G},w}.$ Let $z_I$ be as above. 
Note that for $\zeta \in \mathbb{T}$, there is an absolute constant $a$ so that the following holds. 
\begin{equation} \label{Poisson-trivial}
P_{z_I}(\zeta) \geq a\, \frac{1_I}{|I|}.
\end{equation}

Using Hölder's inequality and \eqref{Poisson-trivial}, we have, 
\begin{align} 
 &\inf_{c}\frac{1}{w(I)}\int_I|\varphi-c|w\:dm\leq \frac{1}{w(I)}\int_I|\varphi-\varphi(z_I)|w\:dm \notag \\
&= \frac{1}{w(I)}\int_I |\varphi -\varphi(z_I)|w\:dm =\frac{1}{w(I)} \frac{|I|}{|I|} \int_I |\varphi -\varphi(z_I)|w\:dm   \notag \\ &=\frac{1}{w_I} \frac{1}{|I|} \int_I |\varphi -\varphi(z_I)|w\:dm 
\leq a^{-1}\frac{1}{w_I} \int_{\mathbb{T}}|\varphi-\varphi(z_I)|w P_{z_I}\:dm \notag \\
&\leq a^{-1} \frac{1}{w_I} \left(\int_\mathbb{T}|\varphi-\varphi(z_I)|^2 w P_{z_I}\:dm\right)^{\frac{1}{2}}\left(\int_{\mathbb{T}}P_{z_I}w\:dm\right)^{\frac{1}{2}} \notag \\
&= a^{-1}\frac{w(z_I)^{\frac{1}{2}}}{w_I} \left(\int_\mathbb{T}|\varphi-\varphi(z_I)|^2P_{z_I}w\:dm\right)^{\frac{1}{2}}. \label{garsia-step}
\end{align}

We now estimate $\frac{w(z_I)^{\frac{1}{2}}}{w_I}$. We do so by first estimating $\frac{w(z_I)}{w_I}$. By Cauchy's inequality, $\frac{1}{w_I}\leq v_I,$ and by \eqref{Poisson-trivial}, $v_I=\frac{1}{|I|}\int_Iw^{-1}\:dm\leq a^{-1}\int_{\mathbb{T}}w^{-1}P_{z_I}\:dm.$ Combining these inequalities we have that 
$$\frac{w(z_I)}{w_I}\leq w(z_I)v_I\leq a^{-1}\int_{\mathbb{T}}wP_{z_I}\:dm\int_\mathbb{T}w^{-1}P_{z_I}\:dm.$$
Finally, since $w \in A_2,$ by Lemma \ref{L:Poisson-Ap}, 
$$\int_{\mathbb{T}}wP_{z_I}\:dm\int_\mathbb{T}w^{-1}P_{z_I}\:dm\leq PA_2(w),$$ and so
\begin{equation}\label{E:interval-and-Poisson}
\frac{w(z_I)}{w_I}\leq a^{-1}PA_2(w).
\end{equation}
From this estimate we have,
\begin{equation} \label{poisson-w-condition}
\frac{(w(z_I))^{\frac{1}{2}}}{w_I}\leq a^{-1}PA_2(w)\frac{1}{(w(z_I))^{\frac{1}{2}}}.
\end{equation}

Returning to \eqref{garsia-step} and using the inequality above, we have
 
\begin{align} \label{garsia-step-ii}
 &\inf_{c}\frac{1}{w(I)}\int_I|\varphi-c|w\:dm\leq a^{-1}\frac{w(z_I)^{\frac{1}{2}}}{w_I} \left(\int_\mathbb{T}|\varphi-\varphi(z_I)|^2P_{z_I}w\:dm\right)^{\frac{1}{2}} \notag \\
 &\leq a^{-2} PA_2(w) \left(\frac{1}{w(z_I)}\int_\mathbb{T}|\varphi-\varphi(z_I)|^2P_{z_I}w\:dm\right)^{\frac{1}{2}}.
\end{align}

By Lemma \ref{L:JN-w-c} with $p=1$, $$||\varphi||_{BMO, w} \approx  \underset{I\subset \mathbb{T}}{\sup} \ \underset{c}{\inf}\left(\frac{1}{w(I)}\int_I |\varphi -c|w\:dm\right),$$ where the implied constants depend on $w.$ Taking the supremum of the left hand side of \eqref{garsia-step-ii}, we conclude that  $$\|\varphi\|_{BMO, w}\lesssim \|\varphi\|_{\mathcal{G}, w}$$ 

We now prove that for each $\varphi \in BMO_{w}$ 
$$ \|\varphi\|_{\mathcal{G}, w} \lesssim \|\varphi\|_{BMO, w},$$ where the implied constants depend on $w$. We do so by following Nikolskii's proof of Lemma \ref{L:Poisson-Ap}, found in \cite{Nik 2002}.
We separately consider intervals with arc length less than $\frac{3}{4}$ and ones with length greater than or equal to $\frac{3}{4}.$ 
Suppose $|I|\geq \frac{3}{4}.$

Using the triangle inequality, we have 
\begin{align} \label{E:triangle-weighted-large}
&\frac{1}{w(z_I)}\int_{\mathbb{T}}|\varphi-\varphi(z_I)|^2P_{z_I}w\:dm \notag \\ &\leq \frac{2}{w(z_I)}\int_{\mathbb{T}}|\varphi-c_{\mathbb{T}}|^2P_{z_I}w\:dm+2|\varphi(z_I)-c_{\mathbb{T}}|^2.
\end{align}
We estimate the two pieces above separately. 

We use the standard estimate of the Poisson kernel, $P_z \leq \frac{(1+|z|)}{1-|z|} \leq \frac {2}{|I|}$. 

Since $P_{z_I}\geq a\frac{1_I}{|I|},$ it follows that $\int_{\mathbb{T}} w P_{z_I}\:dm\geq \frac{a}{|I|} \int_I w\:dm. $ Hence, $\frac{1}{\int_{\mathbb{T}} w P_{z_I}\:dm}\leq \frac{1}{a}\frac{|I|}{w(I)}.$ By combining these inequalities, we have

\begin{align} \label{E:T-poisson}
   \frac{1}{w(z_I)}\int_{\mathbb{T}}|\varphi-c_{\mathbb{T}}|^2P_{z_{I}}w\:dm 
   &\leq \frac{1}{a}\frac{|I|}{\int_{I}w\:dm}\frac{2}{|I|}\int_{\mathbb{T}}|\varphi-c_{\mathbb{T}}|^2w\:dm \notag \\
   &= \frac{2}{a}\frac{1}{\int_Iw\:dm}\int_{\mathbb{T}}|\varphi-c_{\mathbb{T}}|^2w\:dm. 
\end{align}
Taking reciprocals of the inequality given in Lemma \ref{L:Inverse-cond-for-Ap} with $p=2$, with the interval from the lemma taken to be $\mathbb{T}$, and its subset taken to be $I,$ we have
$$\frac{\int_{\mathbb{T}}w\:dm}{\int_{I}w\:dm}\leq [w]_{A_2}\frac{|\mathbb{T}|^2}{|I|^2},$$ and so
$$\frac{1}{\int_Iw\:dm}\leq [w]_{A_2}\frac{1}{|I|^2}\frac{1}{\int_{\mathbb{T}}w\:dm}\leq [w]_{A_2}\frac{16}{9}\frac{1}{\int_{\mathbb{T}}w\:dm},$$
where the last inequality follows from the fact that $|I|\geq \frac{3}{4}.$
Returning to \eqref{E:T-poisson}, and using Lemma \ref{L:JN-weighted}, we have

\begin{align} \label{E:Garsia-weighted-BMO-large}
    \frac{1}{w(z_I)}\int_{\mathbb{T}}|\varphi-c_{\mathbb{T}}|^2P_{z_{I}}w\:dm
    &\leq \frac{2}{a}\frac{1}{\int_{I}w\:dm}\int_{\mathbb{T}}|\varphi-c_{\mathbb{T}}|^2w\:dm \notag \\
   &\leq \frac{32[w]_{A_2}}{9a}\frac{1}{\int_\mathbb{T} w\:dm}\int_{\mathbb{T}}|\varphi-c_{\mathbb{T}}|^2w\:dm  \notag \\
   &\leq \frac{32[w]_{A_2}}{9a}C'_2(w) \|\varphi\|^2_{BMO,w},
\end{align}
with $C'_2(w)$ as in Lemma \ref{L:JN-weighted}.

We now estimate $|\varphi(z_I)-c_{\mathbb{T}}|^2$. Using Hölder's inequality and the estimate \eqref{E:Garsia-weighted-BMO-large}, we have
\begin{align}
    &|\varphi(z_I)-c_{\mathbb{T}}|^2\leq \left(\int_{\mathbb{T}}|\varphi-c_{\mathbb{T}}|P_{z_I}\:dm \right)^2\leq \int_{\mathbb{T}}|\varphi-c_{\mathbb{T}}|^2wP_{z_I}\:dm\int_{\mathbb{T}}w^{-1}P_{z_I}\:dm \notag \\
    &\leq \frac{32[w]_{A_2}}{9a}\int_{\mathbb{T}}wP_{z_I}\:dm\int_{\mathbb{T}}w^{-1}P_{z_I}\:dm \notag \\ &\leq C'_2(w)~\|\varphi\|^2_{BMO,w}\leq C'_2(w)~\frac{32[w]_{A_2}}{9a} ~PA_2(w)~\|\varphi\|^2_{BMO,w}
\end{align}
where the last inequality follows by Lemma \ref{L:Poisson-Ap}.
Returning to \eqref{E:triangle-weighted-large} and combining the estimates above, we have 
\begin{align} 
&\frac{1}{w(z_I)}\int_{\mathbb{T}}|\varphi-\varphi(z_I)|^2P_{z_I}w\:dm \notag \\ &\leq \frac{2}{w(z_I)}\int_{\mathbb{T}}|\varphi-c_{\mathbb{T}}|^2P_{z_I}w\:dm+2|\varphi(z_I)-c_{\mathbb{T}}|^2 \notag \\
&\leq C'_2(w)\frac{64[w]_{A_2}}{9a}(1+PA_2(w))\|\varphi\|^2_{BMO,w}.
\end{align}

We now consider the case $|I|\leq \frac{3}{4}$. Here, we use the estimates from Nikolskii's proof of Lemma \ref{L:Poisson-Ap}, found in \cite{Nik 2002}.
Let $I_k=2^kI/2^{k-1}I$, for $k\geq 1$ and $I_0=I,$ and let $r=|z_I|.$ For $e^{it}\in I_k$, $k \geq 1$
\[|r-e^{it}|=(1-r)^2+4r\sin^2\left(\frac{t}{2}\right)\geq (1-r)^2+16r\frac{2^{k-2}}{\pi^2}\geq 2^{2(k-2)}(1-r)^2,\] 
and for $k=0$, $|r-e^{it}|\geq (1-r)^2.$

By the triangle inequality and Remark \ref{R:induction-BMO}, we have that 
\[|\varphi-c_I|^2\leq 2|\varphi-c_{2^kI}|^2+2|c_{2^kI}-c_I|^2\leq 2|\varphi-c_{2^kI}|^2+2k^24[w]_{A_2}^2\|\varphi\|^2_{BMO, w}.\]

So,  
\begin{align}
    &\int_{\mathbb{T}} |\varphi -c_{I}|^2P_{z_{I}}w\:dm = 
\sum_{k\geq 0}\int_{I_k}|\varphi-c_{I}|^2P_{z_{I}}w\:dm \notag 
\\
&\leq \frac{2}{|I|} \int_{I}|\varphi -c_I|^2w\:dm +2^5 \sum_{k\geq1}\frac{2^{-k}}{|2^kI|}\int_{2^kI}|\varphi-c_I|^2w\:dm. \notag \\
&\leq \frac{2}{|I|}\int_I|\varphi-c_I|^2w\:dm \notag \\ 
&+2^6\sum_{k\geq1}\frac{2^{-k}}{|2^kI|}\int_{2^kI}|\varphi-c_{2^kI}|^2w\:dm\notag \\
&+2^6\sum_{k\geq1}\frac{2^{-k}}{|2^kI|}\int_{2^kI}|c_{2^kI}-c_I|^2w\:dm.\notag
\end{align}
We estimate the three terms in the above sum separately.
Let $$A:=\frac{2}{|I|}\int_I|\varphi-c_I|^2w\:dm,$$
$$B:=2^6\sum_{k\geq1}\frac{2^{-k}}{|2^kI|}\int_{2^kI}|\varphi-c_{2^kI}|^2w\:dm,$$
and 
$$C:=2^6\sum_{k\geq1}\frac{2^{-k}}{|2^kI|}\int_{2^kI}|c_{2^kI}-c_I|^2w\:dm.$$

First we estimate $A.$ 
By Lemma \ref{L:JN-weighted}, we have
\begin{align*}
A=2\frac{w(I)}{|I|}\frac{1}{w(I)}\int_I|\varphi-c_I|^2w\:dm 
&\leq2(C'_1(w))^2\frac{w(I)}{|I|}\|\varphi\|^2_{BMO, w}\\ &\lesssim \frac{w(I)}{|I|}\|\varphi\|^2_{BMO, w},
\end{align*}
with $C'_1(w)$ as in Lemma \ref{L:JN-weighted}. 
Thus, $A \lesssim \frac{w(I)}{|I|}\|\varphi\|^2_{BMO, w}, $ where the implied constant depends on $w.$

We estimate the second term now. First observe that, by Lemma \ref{L:Inverse-cond-for-Ap} and the discussion that follows, we have

\[\left(\frac{|I|}{|2^kI|}\right)^\eta \leq [w]_{A_\eta}\frac{w(I)}{w(2^kI)}, \qquad \text{for  } \eta=2-\epsilon,\: \epsilon >0.\]
Thus, 
\begin{align} \label{eta-condition}
  \frac{w(2^kI)}{|2^kI|}\leq [w]_{A_\eta} \frac{1}{2^{k(1-\eta)}} \frac{w(I)}{|I|} 
\end{align}
Using \eqref{eta-condition} and Lemma \ref{L:JN-weighted}, 
\begin{align}
    B&=2^6\sum_{k\geq1}\frac{2^{-k}}{|2^kI|}\int_{2^kI}|\varphi-c_{2^kI}|^2w\:dm \notag \\ &= 2^6 \sum_{k\geq 1} 2^{-k}  \frac{w(2^kI)}{|2^kI|}\frac{1}{w(2^kI)}\int_{2^kI}|\varphi-c_{2^kI}|^2w\:dm\notag \\ 
    &\leq 2^6 \sum_{k\geq 1} 2^{-k}[w]_{A_\eta}\frac{1}{2^{k(1-\eta)}}\frac{w(I)}{|I|}C'_1(w)\|\varphi\|^2_{BMO, w} \notag \\
    &= 2^6 [w]_{A_\eta}C'_1(w)\sum_{k\geq1}\left(\frac{1}{2^k}\right)^{2-\eta}\frac{w(I)}{|I|}\|\varphi\|^2_{BMO, w}.
\end{align}
Since $\eta = 2-\epsilon <2,$
    it follows that  
    $\lim_{k \to \infty} \left( \frac{\frac{1}{2^{k+1}}}{\frac{1}{2^k}}\right) ^{2-\eta}= 
    \left(\frac{1}{2}\right)^{2-\eta} < 1$.
    So, the sum above converges and 
    \[B \lesssim \frac{w(I)}{|I|}\|\varphi\|^2_{BMO, w},\]
    where the implied constant depends on $w.$

Finally, we have similarly to above, 
\begin{align}
C&=2^6\sum_{k\geq1}\frac{2^{-k}}{|2^kI|}\int_{2^kI}|c_{2^kI}-c_I|^2w\:dm \notag \\
&\leq 2^6 \sum_{k\geq 1} \frac{2^{-k}}{|2^kI|} k^2 c(w)^2 \|\varphi\|^2_{BMO, w} w(2^kI) \notag \\
&= 2^6 c(w)^2 \|\varphi\|^2_{BMO, w} \sum_{k\geq 1} \frac{k^2}{2^k}\frac{2^kI}{|2^kI|}w(2^kI) \notag \\
&\leq 2^6 c(w)^2 \|\varphi\|^2_{BMO, w} \sum_{k\geq1}\frac{k^2}{2^k} [w]_{A_\eta} \frac{1}{2^{k(1-\eta)}}\frac{w(I)}{|I|} \notag \\
&\leq 2^6 c(w)^2 [w]_{A_\eta} \sum_{k\geq1} \frac{k^2}{2^{k(2-\eta)}} \frac{w(I)}{|I|}  \|\varphi\|^2_{BMO, w} \notag.
    \end{align}
Thus, $C\lesssim \frac{w(I)}{|I|}  \|\varphi\|^2_{BMO, w},$ where the implied constant depends on $w.$

By Cauchy-Schwartz, $\int_{\mathbb{T}} w^{-1} P_{z_I}\:dm \int_{\mathbb{T}} w P_{z_I} \:dm \geq 1.$ Using the estimate of the Poisson kernel \eqref{Poisson-trivial},
\begin{align}
\frac{1}{w(z_I)} \frac{w(I)}{|I|} &\leq \int_{\mathbb{T}} w^{-1}P_{z_I}\:dm \frac{w(I)}{|I|} \notag \\
&= \int_{\mathbb{T}} w^{-1} P_{z_I}\:dm \frac{1}{|I|} \int_{I} w \:dm \notag \\
&\leq \int_{\mathbb{T}} w^{-1} P_{z_I} \:dm \frac{1}{a} \int_{\mathbb{T}} w P_{z_I}\:dm \notag \\
&\leq \frac{1}{a} PA_2(w).  \notag 
\end{align}

Therefore, 
\begin{equation} \label{1/w estimate}
\frac{1}{w(z_I)}\leq  \frac{1}{a}PA_2(w) \frac{|I|}{w(I)}.
\end{equation}

Combining the estimates for $A, B$ and $C$, and \eqref{1/w estimate}, we get

\begin{equation} \label{E:A-Garsia-BMO,w}
\frac{1}{w(z_I)}\int_{\mathbb{T}} |\varphi-c_I|^2P_{z_I}w\:dm \lesssim \|\varphi\|^2_{BMO, w}. 
\end{equation}

We estimate the square of the difference between Poisson and weighted averages now. We have
\begin{align*}
    |\varphi-c_I|^2\leq \left(\int_{\mathbb{T}}|\varphi-c_I|P_{z_I}\:dm\right)^2&=\left(\int_{\mathbb{T}}|\varphi-c_I|P_{z_I}w^{\frac{1}{2}}w^{-\frac{1}{2}}\:dm\right)^2 \\
    &\leq \int_{\mathbb{T}}|\varphi-c_I|^2P_{z_I}w\:dm \int_{\mathbb{T}}w^{-1}P_{z_I}\:dm \\
    &\lesssim \|\varphi\|^2_{BMO,w}\int_{\mathbb{T}}wP_{z_I}\:dm\int_{\mathbb{T}}w^{-1}P_{z_I}\:dm \\ 
    &\lesssim \|\varphi\|^2_{BMO,w}PA_2(w) \lesssim  \|\varphi\|^2_{BMO,w}
\end{align*}
where the second inequality follows from Hölder's inequality, the third by estimate \eqref{E:A-Garsia-BMO,w} and the final one by Lemma \ref{L:Poisson-Ap}.

Finally, we combine the estimates above to bound the Garsia norm of $\varphi,$ using the triangle inequality. So,

\begin{align}
&\frac{1}{w(z_I)}\int_{\mathbb{T}}|\varphi-\varphi(z_I)|^2P_{z_I}w\:dm \notag   \\  &\leq \frac{2}{w(z_I)}\int_{\mathbb{T}}|\varphi-c_I|^2P_{z_I}w\:dm + \frac{2}{w(z_I)}\int_{\mathbb{T}}|\varphi(z_I)-c_I|^2P_{z_I}w\:dm  \notag \\
&= \frac{2}{w(z_I)}\int_{\mathbb{T}}|\varphi-c_I|^2P_{z_I}w\:dm + 2|\varphi(z_I)-c_I|^2 \notag \\
&\lesssim (1+PA_2(w))\|\varphi\|^2_{BMO,w}\notag \\
&\lesssim \|\varphi\|^2_{BMO,w}. \notag
\end{align}

We combine the estimates for intervals $|I|\geq \frac{3}{4}$ and $|I|<\frac{3}{4}$ to take the supremum over all the intervals and get
$$\|\varphi\|_{\mathcal{G},w}\lesssim \|\varphi\|_{BMO,w},$$
proving the theorem.

\section{Proof of the Theorem \ref{T:main}}

Having proved Theorem \ref{T:Garsia-weighted}, we now turn to the proof of Theorem \ref{T:main}, in the one-weight case. The proof is similar to Treil's proof of Theorem \ref{T:Bonsall}, but we provide the details for completeness.

\begin{proof}
Let $\varphi$ be the anti-analytic symbol of $H_{\varphi}.$ By Theorem \ref{T:Garsia-weighted}, there exists a constant $C_2(w)$ such that  $\|\varphi\|_{\mathcal{G}}\leq C_2(w) \|\varphi\|_{\mathcal{G},w}.$ 
We may assume, by homogeneity, that $\sup_{z\in \mathbb{D}} \frac{\|H_{\varphi} k_z\|^2_{L^2(w)}}{\|k_z\|^2_{L^2(w)}} = \|\varphi\|^2_{\mathcal{G},w}  \leq \frac{1}{(C_2(w))^2}.$ 
Then, by Theorem \ref{T:Garsia-weighted} and the fact that $\sup_{z\in\mathbb{D}}\frac{\|H_{\varphi} k_{z}\|^2_{L^2(w)}}{\|k_z\|^2_{L^2(w)}}=\|\varphi\|^2_{\mathcal{G}, w},$ it follows that  \begin{align}\|\varphi\|^2_{\mathcal{G}}=\sup_{z\in\mathbb{D}}\int_{\mathbb{T}}|\varphi-\varphi(z)|^2P_z\:dm &\leq (C_2(w))^2      \sup_{z\in\mathbb{D}}\frac{1}{w(z)}\int_{\mathbb{T}}|\varphi-\varphi(z)|^2P_zw\:dm \notag \\ 
&= (C_2(w))^2 \sup_{z\in\mathbb{D}}\frac{\|H_{\varphi} k_{z}\|^2_{L^2(w)}}{\|k_z\|^2_{L^2(w)}} \leq 1. \notag
\end{align}
Rewriting the definition of the Garsia norm, we have
\begin{equation} \label{the-main-harmonic-function} 
\sup_{z \in \mathcal{D}} (|\varphi|^2(z)-|\varphi(z)|^2) \leq 1.
\end{equation}

To estimate the norm of the Hankel operator $H_{\varphi}$ we estimate  
$$(H_{\varphi} f, h)_{L^2(\mathbb{T})} \qquad \text{where } f\in H^2(w),h \in L^2({w^{-1}}).$$ 
By density of analytic polynomials in $H^2(w)$ and density of trigonometric polynomials in $L^2(w^{-1}),$ as well as the definition of the Hankel operator, we may assume that $f$ is an analytic and $h$ trigonometric polynomials. Note that by assumption, $\overline{\varphi}$ is analytic function in $\mathbb{D}$. Then $(H_{\varphi}f,h)_{L^2(\mathbb{T})}=(\varphi f, \overline{g}),$ where $\overline{g}=(I-P)(h).$ Since $h$ is a trigonometric polynomial, $g$ is an analytic polynomial. We may extend this assumption to a disc larger than $\mathbb{D}$ by instead finding the norm of the operator $H_{\varphi_r}$, where $\varphi_r(z)= \varphi (rz) $ for $r \in (0,1),~z\in \mathbb{D}$.  Then, as in \cite{Tre},  we have that
$\overline{\partial} (\varphi f \: g )= f(\overline{\partial}\varphi)g. $
Thus,
$$\Delta(\varphi f\: g )=4(\,\overline{\partial}\varphi\,f'  g+\overline{\partial}\varphi\,fg').$$
Since $(\varphi f \, g)(0)=0,$ by Lemma \ref{Green's Theorem} (Green's theorem) we have that
\begin{align}
(H_{\varphi} f, \overline{g})_{L^2(\mathbb{T})} &= \int_\mathbb{T} \varphi f\: g \:dm = 4\int_{\mathbb{D}}(\overline{\partial}\varphi\:f'  g+\overline{\partial}\varphi\:fg')  \ln \frac{1}{|z|}\:dA(z) \notag \\
&= 4\int_{\mathbb{D}}(\overline{\partial}\varphi\:f'  g+\overline{\partial}\varphi\:fg') w(z)^\frac{1}{2} w(z)^{-\frac{1}{2}}  \ln\frac{1}{|z|}\:dA(z). \label{main-estimate}
\end{align}
We let $dv(z) = \ln \frac{1}{|z|}\:dA(z).$
Using Cauchy-Schwartz's inequality, we have that
\begin{align*}
&\left|4\int_{\mathbb{D}}\overline{\partial}\varphi(z) f(z)\: g'(z)\: w(z)^\frac{1}{2} w(z)^{-\frac{1}{2}}\: dv \right| \\
&\leq \left(4\int_{\mathbb{D}}|\:\overline{\partial}\varphi(z)|^2 |f(z)|^2 w(z)\: dv \right)^{\frac{1}{2}}\left(4\int_{\mathbb{D}} |g'(z)|^2 w^{-1}(z) \: dv \right)^{\frac{1}{2}}. 
\end{align*}
By Lemma \ref{L:weighted-carleson},
\begin{align}
&4\int_{\mathbb{D}}|g'(z)|w^{-1}(z) \: dv \leq 4c(w^{-1})^2 \int_{\mathbb{{T}}}|g-\hat{g}(0)|^2w^{-1}\:dm \notag \\
&\leq 4c(w^{-1})^2 \left(\int_{\mathbb{T}} |g|^2w^{-1}\:dm + |\hat{g}(0)|^2\int_{\mathbb{T}} w^{-1}\:dm\right). \notag
\end{align}
Using Hölder's inequality, we have
\begin{align}
    |\hat{g}(0)|^2 &= \left|\int_{\mathbb{T}}g\:dm\right|^2 
    = \left(\int_{\mathbb{T}}|g|^2w^{-1}\:dm\right)\left(\int_{\mathbb{T}}w\:dm\right). \notag 
\end{align}
Thus, \begin{align*}
4\int_{\mathbb{D}}|g'(z)|w^{-1}\: dv \leq 16c(w^{-1})^2[w]_{A_2}\int_{\mathbb{T}}|g|^2w^{-1}\:dm.
\end{align*}
We let $u(z)=: 1+ |\varphi(z)|^2-|\varphi|^2(z).$ Note that $u(z)$ is subharmonic. By \eqref{the-main-harmonic-function} it follows that $0 \leq u(z) \leq 1.$ Thus, Lemma \ref{L:Uchiyama}, now implies that $d\tau(z):=\Delta u \ln \frac{1}{|z|} \:dA(z)=4|\:\overline{\partial}\varphi(z)|^2\ln\frac{1}{|z|}\:dA(z) $ is a Carleson measure. Since $d\tau$ is Carleson, Lemma \ref{L:weighted-carleson} then implies that
\begin{align*}
    4\int_{\mathbb{D}}\left|\overline{\partial}\varphi\right(z)|^2|f(z)|^2 w(z) \: dv \leq C_{\tau}(w) \|f\|^2_{L^2(w)}.
\end{align*}
Together, these estimates give 
\begin{align*}\left|4\int_{\mathbb{D}}\overline{\partial}\varphi(z)f(z) g'(z)\: w^\frac{1}{2}(z) w^{-\frac{1}{2}}  \: dv\right| \leq \sqrt{C_{\tau}(w)}\,c(w^{-1})\,\|f\|_{L^2(w)}\|g\|_{L^2(w^{-1})}.
\end{align*}
We use Cauchy-Schwartz to estimate the second term in \eqref{main-estimate}. So, we have
\begin{align}
&\left|4\int_{\mathbb{D}}\overline{\partial}\varphi(z)\:f'(z)g(z) w^\frac{1}{2}(z) w^{-\frac{1}{2}} (z)  \: dv\right| \notag \\ &\leq \left(4\int_{\mathbb{D}}|\overline{\partial}\varphi(z)|^2 |g(z)|^2 w^{-1}(z) \: dv \right)^{\frac{1}{2}}\left(4\int_{\mathbb{D}} |f'(z)|^2 w(z) \: dv\right)^{\frac{1}{2}}. \notag
\end{align}
Changing the role of $f$ and $g$ above, we get 
$$
\left|4\int_{\mathbb{D}}\overline{\partial}\varphi(z)\:f'(z)g(z) w^\frac{1}{2}(z) w^{-\frac{1}{2}} (z)  \: dv \right| \leq \sqrt{C_{\tau}(w)}\, c(w)\, \|f\|_{L^2(w)}\|g\|_{L^2(w^{-1})}.
$$
Tracking constants in Nikolskii's proof of Lemma \ref{L:weighted-carleson} (property f)ii) on Page 111 in \cite{Nik 2002}), we have that $c(w)=c(w^{-1})$. Finally, 
$$|(H_{\varphi} f, \overline{g})_{L^2(\mathbb{T})}|\leq \sqrt{C_{\tau}(w)}\, c(w)\, \|f\|_{L^2(w)}\|\overline{g}\|_{L^2(w^{-1})}.$$

Finally, since $I-\mathbb{P}$ is bounded on $L^2(w^{-1}),$ for $w^{-1}\in A_2,$ we have that $\|\overline{g}\|_{L^2(w^{-1})}=\|(I-\mathbb{P})(h)\|_{L^2(w^{-1})}\lesssim \|h\|_{L^2(w^{-1})}.$

Therefore, $H_{\varphi}$ is bounded. 

By homogeneity, we can conclude that $\|H_{\varphi}\|_{L^2(w)} \lesssim \|\varphi\|_{\mathcal{G},w}$, where the implied constants depend on the $A_2$ characteristic of $w.$

\end{proof}

\begin{section}{Two-weight Carleson theorem and the Proof of Theorem \ref{T:main}}
In this section, we prove the two-weight Carleson theorem as well as Theorem \ref{T:main} for two weights. The proof of the theorem of the Carleson's embedding theorem for two weights uses similar ideas to the proof of the unweighted Carleson embedding theorem, found in Garnett \cite{Gar}.

\begin{subsection}{Proof of two-weight Carleson embedding theorem} 
\begin{proof}
Note that $i)$ implies $ii),$ since for each $z\in \mathbb{D},$ $k_z$ is bounded and so belongs to $L^2(\mu).$

We now prove that $ii)$ implies $i).$

By assumption,
$$\int_{\mathbb{D}}P_z(\zeta) (\mu(\zeta))^2 \frac{1}{(\mu(\zeta))^2} \lambda(\zeta)\:d\tau(\zeta) \leq B \int_{\mathbb{T}}P_z(\xi)\mu(\xi)dm(\xi).$$
Since $\mu \in A_2,$ by Lemma \ref{L:Poisson-Ap}, it follows that $\mu^{-1}(\zeta)\leq [\mu]_{A_2} \frac{1}{\mu(\zeta)}.$ Thus,
$$\int_{\mathbb{D}}P_z(\zeta)(\mu(\zeta))^2(\mu^{-1}(\zeta))^2 \lambda(\zeta) \: d\tau(\zeta)\leq  B\, \lbrack \mu \rbrack_{A_2} \int_{\mathbb{T}}P_z(\xi) \mu(\xi) \: dm(\xi)$$ 
Since $\mu \in A_2,$ by Lemma \ref{L:One-weight-Treil-and-Volberg}, it follows that 
$$\int_{\mathbb{D}}|f(\zeta)|^2 (\mu^{-1}(\zeta))^2\lambda(\zeta) \:d\tau(\zeta) \lesssim B\, [\mu]_{A_2} \int_{\mathbb{T}}|f(\xi)|^2\mu^{-1}(\xi)\: dm(\xi),$$ for all $f\in L^2(\mathbb{T},\mu^{-1}).$
In particular, for each $z\in \mathbb{D},$ $k_z \in L^2(\mathbb{T},\lambda^{-1}),$ since for each $z\in \mathbb{D},$ $k_z$ is bounded. Therefore,
$$\int_{\mathbb{D}}P_z(\zeta)(\mu^{-1}(\zeta))^2 \lambda (\zeta) \: d\tau(\zeta) \lesssim B\, [\mu]_{A_2}\int_{\mathbb{T}}P_z(\xi)\mu^{-1}(\xi) \: dm.$$
Applying Lemma \ref{L:One-weight-Treil-and-Volberg} again with $\mu^{-1}$ we get,
$$\int_{\mathbb{D}}|f(\zeta)|^2\lambda(\zeta)\:d\tau(\zeta)\lesssim B\, [\mu]_{A_2}\int_{\mathbb{T}}|f(\xi)|^2\mu(\xi)\:dm(\xi),$$
for all $f\in L^2(\mathbb{T},\mu)$

To complete the proof of the above theorem, we show that $ii)$ and $iii)$ are equivalent.

Assume that $ii)$ holds. Suppose $S=\{re^{i\theta}: 1-h\leq r\leq 1, |\theta-\theta_0|<h\}.$ Taking $z_0=0,$ 
we have that $\frac{1}{\mu(\mathbb{T})}(\lambda \: d\tau)(\mathbb{D})\leq B. $ Suppose $h \geq \frac{1}{4},$ and let $S$ be a sector as above. Then
\begin{align*}
\frac{1}{\mu(I_S)}(\lambda \: d\tau)(S)&\leq \frac{1}{\mu(I_S)} (\lambda \: d\tau)(\mathbb{D}) \\ &\leq B \, \frac{\mu(\mathbb{T})}{\mu(I_S)} \leq B\, C(\mu) \left(\frac{|\mathbb{T}|}{|I_S|}\right)^2 \leq B\, C(\mu) \left(\frac{1}{h}\right)^2 \lesssim B\, C(\mu), 
\end{align*}
since $h \geq \frac{1}{4}.$
Assume $h< \frac{1}{4},$ let $S$ be a sector of the disc as above.  Let $z_0 = (1-\frac{1}{2}h)e^{i\theta_0}.$ For 
$z\in S,$ we have 
$$\frac{1-|z_0|^2}{|1-\overline{z_0}z|^2}\geq \frac{C}{1-|z_0|^2}.$$
Therefore, 
\begin{align}
\frac{1}{\mu(I_S)}(\lambda \: d\tau)(S) 
&\lesssim \frac{(1-|z_0|)}{\mu(I_S)} \int_{\mathbb{D}}\frac{1-|z_0|^2}{|1-\overline{z_0}z|^2} \lambda \: d\tau \lesssim B\, \frac{(1-|z_0|)}{\mu(I_S)} \mu (z_0) \notag \\ &= B\, \frac{1}{2}\frac{|I_S|}{\mu(I_S)} \mu(z_0) \lesssim B \, C(\mu). \notag
\end{align}
The interval $I_{z_0},$ defined as in the proof of Theorem \ref{T:Garsia-weighted}, and the interval $I_S$ have the same center while the interval $I_{z_0}$ is half the length of $I_S.$ Thus we can use the equation (\ref{E:interval-and-Poisson}) to conclude that the last inequality holds.
Therefore $ii)$ implies $iii)$

Suppose now that $iii)$ holds. Let $z_0=re^{i\theta_0} \in \mathbb{D}.$ Suppose $|z_0|< \frac{3}{4}.$ Let $I_0$ be an interval centered at $e^{i\theta}$ of length $1-|z_0|.$ Since $P_{z_0}\geq a \frac{1_{I_0}}{|I_0|},$ where $a$ is an absolute constant, it follows that
$\frac{1}{\mu(z_0)}\lesssim \frac{|I_0|}{\mu(I_0)}.$
So, 
\begin{align*}\frac{1}{\mu(z_0)}\int_{\mathbb{D}}P_{z_0}(z)\lambda \: d\tau \lesssim D \, \frac{\mu(\mathbb{T})}{\mu(z_0)}\lesssim D \, \mu(\mathbb{T})\frac{|I_0|}{\mu(I_0)} &\lesssim D \, C(\mu) \left(\frac{|\mathbb{T}|}{|I_0|}\right)^2 |I_0| \\ &= D\, C(\mu)\frac{1}{|I_0|} \lesssim D \, C(\mu).\end{align*}
Assume $|z_0|> \frac{3}{4}.$ Let 
$$B_n = \{z\in \mathbb{D}: |z-\left(\frac{z_0}{|z_0|}\right)|< 2^n (1-|z_0|)\}.$$
Let $I_{n}$ be an interval centered at $\frac{z_0}{|z_0|}$ of length $2^n(1-|z_0|),$ for $n=1,2,3,\ldots .$
By assumption, 
$$(\lambda \: d\tau)(B_n) \leq D \mu (I_{n}).$$
For $z\in B_1,$
$$\frac{1-|z_0|^2}{|1-\overline{z_0}z|^2}\leq \frac{C}{1-|z_0|},$$
and for $n\geq 2,$ $z\in B_n \setminus B_{n-1},$
$$\frac{1-|z_0|^2}{|1-\overline{z_0}z|^2}\leq \frac{C}{2^{2n}(1-|z_0|)}.$$
Therefore, 
\begin{align}
\int_{\mathbb{D}}\frac{1-|z_0|^2}{|1-\overline{z_0}z|^2} \lambda d\tau &\leq \int_{B_1} \frac{1-|z_0|^2}{|1-\overline{z_0}z|^2} \lambda(z) \: d\tau(z) + \sum_{n=2}^{\infty} \int_{B_n}\frac{1-|z_0|^2}{|1-\overline{z_0}z|^2} \lambda(z) \: d\tau(z) \notag \\
&\lesssim \sum_{n=1}^{\infty} \frac{(\lambda \: d\tau)(B_n)} {2^{2n}(1-|z_0|)} \leq D \sum_{n=1}^{\infty} \frac{\mu(I_{n})}{|I_{n}|}\frac{1}{2^n} \notag \\ &\lesssim  D \: C(\mu)\frac{\mu(I_0)}{|I_0|}\sum_{n=1}^{\infty} 2^{n(\eta-2)},\notag
\end{align}
where the last inequality follows from Lemma \ref{L:Inverse-cond-for-Ap}, since $\mu \in A_{\eta},$ with $\eta<2. $ Since $\eta<2,$ the last sum converges. 
As above $\frac{1}{\mu(z_0)}\lesssim \frac{|I_0|}{\mu(I_0)},$ and therefore 
$$\frac{1}{\mu(z_0)} \int_{\mathbb{D}}\frac{1-|z_0|^2}{|1-\overline{z_0}z|^2} \lambda d\tau < C \: D(\mu),$$
proving the theorem.

\end{proof}

\end{subsection}

\begin{subsection}{Boundedness of Hankel Operators in the two-weight setting}

We first prove that under the assumptions of Theorem \ref{T:main}, Poisson integrals of a certain non-negative measure are uniformly bounded. 

\begin{prop}  \label{P:Poisson-Condition-1-z} 
Let $\mu, \lambda$ be Muckenhoupt weights on $\mathbb{T}.$ 
Suppose 
$$\sup_{z_0\in \mathbb{D}}\frac{1}{\mu(z_0)} \int_{\mathbb{T}}|\varphi-\varphi(z_0)|^2\lambda \: dm < \infty.$$ Then,

\begin{align*}&\sup_{z_0\in \mathbb{D}}\frac{1}{\mu(z_0)} \int_{\mathbb{D}}|\nabla\varphi(z)|^2 \lambda (z)(1-|z|^2)P_{z_0}(z)\:dA(z) \\ &\lesssim \sup_{z_0\in \mathbb{D}}\frac{1}{\mu(z_0)}\int_{\mathbb{T}} |\varphi-\varphi(z_0)|^2\lambda \: dm,\end{align*}

and the measure $(1-|z|^2)\:dA(z)$ is $\lambda-\mu$ Carleson.
\end{prop}

\begin{proof}
    Replacing $f$ by $\frac{\varphi \circ b_{z_0}}{(\mu(z_0))^{\frac{1}{2}}}- \frac{\varphi \circ b_{z_0}(0)}{(\mu(z_0))^{\frac{1}{2}}},$ where $b_{z_0}$ is an automorphism of the circle that interchanges 0 to $z_0, $ in Lemma \ref{L:Green-weighted}, after a change of variables we get
    $$\frac{1}{\mu(z_0)}\int_{\mathbb{T}}|\varphi-\varphi(z_0)|^2 \lambda \: dm \approx \frac{1}{\mu(z_0)} \int_{\mathbb{D}} |\nabla{\varphi(z)}|^2 \lambda(z) \ln \frac{1}{|b_{z_0}(z)|} \: dA(z).$$
Since $1-|z|^2 \lesssim \ln \frac{1}{|z|}$ for all $z\in \mathbb{D},$ we have 
$(1-|z|^2)P_{z_0}(z)=1-|b_{z_0}(z)|^2 \lesssim \ln \frac{1}{|b_{z_0}(z)|,}$ for all $z\in\mathbb{D}.$ Therefore, 
\begin{align}&\frac{1}{\mu(z_0)} \int_{\mathbb{D}}|\nabla\varphi(z)|^2 \lambda (z)(1-|z|^2)P_{z_0}(z)\:dA(z) \notag \\  &\lesssim \frac{1}{\mu(z_0)}\int_{\mathbb{D}} |\nabla{\varphi(z)}|^2 \lambda(z) \ln \frac{1}{|b_{z_0}(z)|} \: dA(z) \notag \\  &\approx \frac{1}{\mu(z_0)}\int_{\mathbb{T}} |\varphi-\varphi(z_0)|^2\lambda \: dm. \notag
\end{align}
By the condition $ii)$ of Theorem \ref{T:Two-weight-Carleson}, it then follows that the measure $(1-|z|^2)\:dA(z)$ is $\lambda-\mu$ Carleson.
\end{proof}

To prove that Hankel operators are bounded under the assumptions of Theorem \ref{T:main}, we first prove that the measure $\ln(\frac{1}{|z|})\:dA(z)$ is two-weight Carleson. We do so using the fact that the measure $(1-|z|^2)\:dA(z)$ is $\lambda-\mu$ Carleson. Our proof uses subharmonicity of the appopriate function to link these two cases, following the ideas in Garnett \cite{Gar}.

\begin{prop} \label{P:Green-Carleson-measure} Let $\varphi$ be an analytic  or anti-analytic function in $L^1(\mathbb{T}).$
Let $$dv_1(z)=|\nabla \varphi(z)|^2 \log \frac{1}{|z|}\: dA(z),$$
$$dv_2(z) = |\nabla \varphi(z)|^2 (1-|z|^2)\: dA(z).$$ 
Let $\mu, \lambda \in A_2.$
Then $dv_1$ is $\lambda -\mu $ Carleson if and only if $dv_2$ is $\lambda -\mu $ Carleson. 
\end{prop}

\begin{proof}
Suppose $\varphi$ is anti-analytic. 
    The inequality $1-|z|^2 \leq 2 \log (\frac{1}{|z|}),$ for all $z\in \mathbb{D}$ shows that $dv_1$ being $\lambda -\mu $ Carleson implies that $dv_2$ is \mbox{$\lambda-\mu$} Carleson. 
For $|z|> \frac{1}{4},$ there exists a constant $C$ such that the reverse inequality \mbox{$\log\frac{1}{|z|}\leq C\, (1-|z|^2)$} holds. 
If we assume that $h\leq \frac{3}{4},$ the inequality above implies that for sectors $S=\{re^{i\theta}: 1-h \leq r <1, |\theta-\theta_0|<h\},$ \mbox{$dv_1(S)\leq C \,dv_2(S).$} Then $(\lambda \: dv_1)(S) \leq C\, (\lambda \: dv_2)(S).$ So the inequality holds for $h\leq \frac{3}{4}.$ We now prove that $(\lambda \: dv_1)(\{z: |z|\leq \frac{1}{4}\}) \leq C\, (\lambda \: dv_2) ({z: |z|\leq \frac{1}{2}}). $ These two estimates prove that $dv_2$ is Carleson. 
To prove the second estimate, we prove that the function $\frac{1}{\lambda^{-1}}|\nabla \varphi|^2$ is subharmonic. We let \mbox{$v(z)=-\log(\lambda^{-1}(z)).$}Since $\varphi$ is antianalytic, $|\nabla \varphi|^2=|\overline{\partial}\varphi|^2, $ and so
 \begin{align}
 &\Delta\left(\frac{1}{\lambda^{-1}(z)} |\overline{\partial}\varphi(z)|^2\right)= \Delta(e^{v(z)}|\overline{\partial}\varphi(z)|^2)) \notag \\ &= e^{v(z)} (\Delta(v(z))|\overline{\partial}\varphi(z)|^2+4|\overline{\partial}(v(z))\overline{\partial}\varphi(z)+\overline{\partial}\overline{\partial}\varphi(z)|^2) \notag \\ &\geq e^{v(z)} (\Delta(v(z))|\overline{\partial}\varphi(z)|^2) \notag \\
 &= \frac{1}{\lambda^{-1}(z)}\left(\frac{|\nabla \lambda(z)|^2}{(\lambda(z))^2}\right) |\overline{\partial}\varphi(z)|^2 \geq 0. \notag
 \end{align}

 Therefore $\frac{1}{\lambda^{-1}(z)}|\overline{\partial}\varphi(z)|^2$ is subharmonic.
 
 Now, 
 \begin{align*}(\lambda \: dv_1)(\{z: |z|\leq \frac{1}{4}\}) &\leq C(\lambda) \sup_{|z|\leq\frac{1}{4}} \frac{1}{\lambda^{-1}(z)}|\nabla \varphi(z)|^2 \\ &\lesssim C(\lambda) \int_{|\zeta|< \frac{1}{2}} |\nabla \varphi(\zeta)|^2 \frac{1}{\lambda^{-1}(\zeta)} (1-|\zeta|)  d\zeta d\eta  \\ 
 &\lesssim C(\lambda) (\lambda \: dv_2)(\{\zeta: |\zeta|\leq \frac{1}{2}\}).
 \end{align*}

 Thus $dv_1$ is $\lambda-\mu$ Carleson. 
 
 Similar proof holds if $\varphi$ analytic. 

\end{proof}

We are now ready to prove Theorem \ref{T:main}. The proof is similar to the one-weight case and so follows the approach of Treil \cite{Tre}. We therefore omit some of the details.

\begin{proof} [Proof of Theorem \ref{T:main}]
  Propositions \ref{P:Poisson-Condition-1-z} proves that under the assumption i), 
  the measure $dv_1(z)=|\nabla \varphi (z)|^2 (1-|z|^2) \:dA(z),$ is $\lambda-\mu$ Carleson. Proposition \ref{P:Green-Carleson-measure} then implies that the measure 
  $|\nabla \varphi (z)|^2 \log\frac{1}{|z|} \: dA(z)$ is $\lambda-\mu$ Carleson. Thus,  for each $f\in L^2(\mathbb{T},\mu \: dm),$ 
    we have  
    \begin{equation} \label{E:Carleson-mu}\int_{\mathbb{D}}|f(z)|^2\lambda(z)\: |\nabla \varphi|^2\ln{\frac{1}{|z|}}\:dA(z) \lesssim \int_{\mathbb{T}}|f(\xi)|^2 \mu(\xi) \: dm(\xi).
    \end{equation}
    Similarly to above, condition $ii)$ implies that for each $g\in L^2(\mathbb{T}, \lambda^{-1}),$
    \begin{equation} \label{E:Carleson-mu-conjugate}
    \int_{\mathbb{D}}|g(z)|^2\mu^{-1}(z)|\nabla \varphi|^2\ln{\frac{1}{|z|}} \: dA(z) \lesssim \int_{\mathbb{T}}|g(\xi)|^2 \lambda^{-1}(\xi) \: dm(\xi).
    \end{equation}

In particular, tracing the dependence on the constants above, we see that (\ref{E:Carleson-mu}) becomes 
\begin{equation} \label{E:Carleson-mu-constants}\int_{\mathbb{D}}|f(z)|^2\lambda(z) |\nabla \varphi|^2\ln{\frac{1}{|z|}}\:dA(z) \lesssim \|\varphi\|^2_{\mathcal{G},\lambda\mu} \int_{\mathbb{T}}|f(\xi)|^2 \mu(\xi) \: dm(\xi)
\end{equation}
and (\ref{E:Carleson-mu-conjugate}) becomes

\begin{equation} \label{E:Carleson-mu-conjugate-constants}
    \int_{\mathbb{D}}|g(z)|^2\mu^{-1}(z)|\nabla \varphi|^2\ln{\frac{1}{|z|}}\: dA(z) \lesssim  \|\varphi\|^2_{\mathcal{G},\mu^{-1},\lambda^{-1}}\int_{\mathbb{T}}|g(\xi)|^2 \lambda^{-1}(\xi) \: dm(\xi).
    \end{equation}

    We use the facts above and duality to estimate $|(H_{\varphi}f, h)_{L^2(\mathbb{T})}|$, for $f \in L^2 (\mathbb{T}, \mu \: dm), h\in L^2(\mathbb{T}, \lambda^{-1} \: dm).$ 
    
    We assume $f$ is an analytic polynomial, and $h$ trigonometric one. As before, we let $\overline{g}=(I-\mathbb{P})(h).$ Then both $f,\overline{g}$ are analytic polynomials. To prove the above inequality we prove that $$|(H_{\varphi}f, \overline{g})_{L^2(\mathbb{T})}| \lesssim \max\{\|\varphi\|_{\mathcal{G},\lambda\mu}, \|\varphi\|_{\mathcal{G},\mu^{-1}\lambda^{-1}}\}\|f\|_{L^2(\mathbb{T}, \mu \: dm)} \|g\|_{L^2(\mathbb{T}, \lambda^{-1}\: dm)}.$$
   The result then follows by density of analytic polynomials in 
    $H^2(\mathbb{T},\lambda^{-1} \: dm)$
    and trigonometric polynomials in $L^2(\mathbb{T}, \mu \: dm),$ as well as the boundedness of the operator $I-\mathbb{P}$ on $L^2(\mathbb{T},\, \lambda^{-1} \: dm),$ for $\lambda^{-1}\in A_2.$ 
Similarly to the one-weight case, we use Green's theorem and the fact that
$\overline{\partial} (\varphi f \: g )= f(\overline{\partial}\varphi)g, $
to estimate the norms on the disk. We apply Cauchy-Schwartz with $\lambda$ instead of $w$ and use (\ref{E:Carleson-mu}) and Lemma \ref{L:Green-weighted} to estimate the first term and (\ref{E:Carleson-mu-conjugate}) and Lemma \ref{L:Green-weighted} to estimate the second term. This gives the desired estimates for $|(H_{\varphi}f,\overline{g})|$ and allows us to conclude that $\|H_{\varphi}\|_{L^2(\mathbb{T}, \lambda)\to L^2(\mathbb{T}, \mu)} \lesssim \max\{\|\varphi\|_{\mathcal{G},\lambda\mu}, \|\varphi\|_{\mathcal{G},\mu^{-1}\lambda^{-1}}\},$ where the implied constants depend on $\mu$ and $\lambda.$

Letting $f\in L^2(\mathbb{T}, \lambda^{-1}), g\in L^2(\mathbb{T}, \mu),$ we can similiarly prove that 
$$\|H_{\varphi}\|_{L^2(\mathbb{T},\, \lambda^{-1})\to L^2(\mathbb{T},\, \mu^{-1})} \lesssim \max\{\|\varphi\|_{\mathcal{G},\lambda\mu}, \|\varphi\|_{\mathcal{G},\mu^{-1}\lambda^{-1}}\}. $$ Thus the norm of the Hankel operators between the weighted spaces is controlled by the two-weight Garsia norm of its symbol.
\end{proof}

\end{subsection}

\end{section}

\begin{section}{Application of Theorem \ref{T:Two-weight-Carleson} to a class of integral operators}

We characterize two-weight boundedness of a certain class of integral operators described by Aleman and Siskakis in \cite{AleSis}. 
We define an integral operator on $H^2(\mu)$ densely on polynomials. 

For $f$ a polynomial in $L^2(\mu)$,
let 
$$T(f)(z)=T_g(f)(z)=\frac{1}{z}\int_{0}^{z}f(\zeta)g'(\zeta)d\zeta,$$
for a function $g$ analytic on $\mathbb{D}.$ We also define $P_g=M_zT_g,$ for $z\in\mathbb{D},$ where $M_z$ is a multiplication operator. 
We estimate the norm of $T_g$ from $H^2(\mathbb{T}, \mu)$ to $H^2(\mathbb{T}, \lambda)$ where $\mu$ and $\lambda$ are $A_2$ weights. Observe that $T_g$ is bounded between these spaces if and only if $P_g$ is bounded. So, we study the boundedness of $P_g$ instead.  
We have the following corollary of Theorem \ref{T:Two-weight-Carleson}. 

\begin{cor}
 Let $\mu,\lambda$ be weights in $A_2.$ Let $g$ be a function in $H^2(\mathbb{T}).$ Let $T_g$ be given by $$T(f)(z)=T_g(f)(z)=\frac{1}{z}\int_{0}^{z}f(\zeta)g'(\zeta)d\zeta,$$ on analytic polynomials $f$.
\begin{itemize}
\item[i)] Suppose that $T_g$ is defined as above and that the function $g$ satisfies:
$$\left(\sup_{z\in\mathbb{D}}\frac{1}{\mu(z)}\int_{\mathbb{T}}|g-g(z)|^2 P_z \: \lambda \:dm\right)^{\frac{1}{2}} < \infty.$$
Then $T_g$ is bounded as an operator between $H^2(\mathbb{T},\mu)$ and $L^2(\mathbb{T},\lambda).$
\item[ii)] 
Suppose that $T_g$ is defined as above, and that the function $g$ satisfies:
$$\left(\sup_{z\in\mathbb{D}}\frac{1}{\lambda^{-1}(z)}\int_{\mathbb{T}}|g-g(z)|^2 P_z \: \mu^{-1}\:dm\right)^{\frac{1}{2}} < \infty,$$
$T_g$ is bounded as an operator between $H^2(\mathbb{T},\lambda^{-1})$ and $L^2(\mathbb{T},\mu^{-1}).$
\end{itemize}
\end{cor}

\begin{proof}
We prove $i).$ Statement $ii)$ follows similarly. 
By the remarks above, we estimate the norm of the operator $P_g$ instead. We use duality, and density of polynomials as in the proof of boundedness of Hankel operators. 
For $g\in H^2$ and $f,h$ polynomials, the Littlewood Paley identity gives 

\begin{align*}
|(P_g f, h)_{L^2(\mathbb{T})}| &= \left|2\int_{\mathbb{D}}(P_gf)'(z)\overline{h'(z)} dv(z) \right| \\ &= \left|2 \int_{\mathbb{D}} f(z)g'(z)\overline{h'(z)} dv(z) \right|,
\end{align*}
where $dv(z)= \log\frac{1}{|z|}\:dA(z).$
Note that by assumption and Theorem \ref{T:Two-weight-Carleson}, we have that $|g'(z)|^2\log\frac{1}{|z|}\: dA(z)$ is $\lambda-\mu$ Carleson. 
Using Cauchy-Schwartz and the fact that $\lambda^{-1}(z)\approx \frac{1}{\lambda(z)},$ for $z\in \mathbb{D},$ and $\lambda \in A_2$ we have
\begin{align*}
    \left|(P_g f, h)_{L^2(\mathbb{T})}\right| &\lesssim \left(\int_{\mathbb{D}}|f(z)|^2\lambda(z)|g'(z)|^2\:dv(z)\right)^{\frac{1}{2}}\left(\int_{\mathbb{D}}|h(z)|^2 \lambda^{-1}(z)\: dv(z)\right)^{\frac{1}{2}} \\
    &\lesssim \left(\int_{\mathbb{T}}|f(\zeta)|^2\mu(\zeta)\: dm(\zeta)\right)^{\frac{1}{2}} \left(\int_{\mathbb{T}}|h(\zeta)|^2\lambda^{-1}(\zeta)\: dm(\zeta)\right)^{\frac{1}{2}},
\end{align*}
where the last inequality follows by Theorem \ref{T:Two-weight-Carleson}, and Lemma \ref{L:Green-weighted}. Thus the operator $P_g$ is bounded between $H_2(\mu)$ and $H_2(\lambda),$ and therefore $T_g$ is bounded between these spaces as well.  

\end{proof}

\end{section}

\subsection*{Acknowledgements}
    The author would like to thank Brett Wick for helpful discussions.

\end{document}